\documentclass{autart_nojtitle}
\usepackage{amssymb}
\usepackage{amsmath}
\usepackage{graphicx}
\usepackage{pict2e} %  for platex
\newcommand{\calN}{\mathcal{N}}

\newcommand{\bfI}{\mathbf{I}}

\newcommand{\bfK}{\mathbf{K}}

\begin{document}
\begin{frontmatter}
\runtitle{Modification of max-separable function for iISS}
\title{Interval Observer of Minimal Error Dynamics\thanksref{footnoteinfo}}
\thanks[footnoteinfo]{
The work was supported in part by JSPS KAKENHI Grant Number 17K06499. 
The material in this paper was not presented at any conference. \\
$^\ast$ Corresponding author. 
Tel.+81-948-29-7717. Fax +81-948-29-7709. 
} 
\author[KIT]{Hiroshi Ito}\ead{hiroshi@ces.kyutech.ac.jp}$^{,\ast}$, 
%\author[KAG]{}\ead{}, 
\address[KIT]{
Department of Intelligent and Control Systems, 
Kyushu Institute of Technology, 
680-4 Kawazu, Iizuka, 820-8502, Japan}%
%\address[KAG]{}%
%\subtitle{}    % 35 words or less. Regulars paper. 
%
\begin{abstract} 
This paper proposes a simple interval observer which can 
generate tighter interval estimates of variables 
in transient states than the standard interval observer. 
The simple nonlinear dynamics shrinks the estimated intervals 
to true state variables at the maximum velocity in the absence of 
disturbances. 
In the presence of bounded disturbances, 
ultimate bounded of the interval estimates are 
given under the standard assumption for interval observer design. 

\end{abstract}
\end{frontmatter}
%% body %%%
%=============================================
%\thispagestyle{plain}  % for test
%\pagestyle{plain}      % for test
%\renewcommand{\thefootnote}{\fnsymbol{footnote}}
%
%\addtolength{\baselineskip}{-.15mm}
%
\section{Introduction}\label{sec:intro}

In contrast to classical observers which estimate unmeasured variables of dynamical systems asymptotically, interval observers provide intervals to which unmeasured variables are guaranteed to belong all the time \cite{GouzeRapaport_EcMo_00,Mazebern,efimrz}. 
The set estimation achieved by interval observers in transient 
phases are beneficial in the presence of non-stationary disturbances. 
Interval observers are found to be useful for state estimation in many applications (see \cite{olijl,GouzeRapaport_EcMo_00,mobego,intobsDAMADICS06,intobsCANAL10,intobsTURBINE14} 
to name a few). Reading \cite{EfiRassiSurv16} may allow one to quickly overview interval observers developed for a variety of systems. This brief paper is not in the direction of expanding system classes to cover. This paper focuses on the tightness of interval estimates. 
As demonstrated in \cite{olijl,mobego}, to tighten interval estimates, one can build a lot of interval observers and take the intersection of all the estimated intervals. Although this approach is actually effective in practice, the total size of multiple interval observes grows rapidly as one wants better estimates. The multiple observers work independently, so that a better estimate at one moment is not exploited to produce estimates at another moment. 
Reinitialization is an idea to ease the disorder of the estimates generated by the individual observes running in parallel \cite{mobego}, although  
the reinitialization mechanism taking place at time intervals chosen somehow is not yet completely efficient for 
the tightening in spite of its complexity. 
Another approach is to search an observer gain that minimizes 
a criterion whose scalar value somehow represents the width of 
all the intervals \cite{EfiRassiSurv16}. 
This parameter tuning restricts the observer dynamics to 
the structure of a single linear observer gain. 
This paper removes this restriction and introduces a simple nonlinear 
mechanism of tightening the interval estimates 
without increasing the dimension of the observer by mixing mechanisms of 
estimate update without reinitialization. 
This paper proves that the simple idea indeed gives an 
interval observer 
with tighter estimates by invoking the monotonicity property of 
error systems.

\noindent
{\it Notation: }
In this paper, $\Rset$ denotes the set of real numbers. 
The set of non-negative real numbers is denoted by 
$\Rset_+$, i.e., $\Rset_+:=[0,\infty)$. 
The symbol $|v|$ denotes the Euclidean norm
for $v\in\Rset^n:=(-\infty,\infty)^n$. 
For $x, y\in \Rset^n$, 
$x \le y$ if $y-x \in \Rset^n_+$. 
We write $x < y$ if $x \le y$  and $x \neq y$. 
The expression $x \ll y$ is used if $y-x$ is in the interior of $\Rset^n_+$. 
A square matrix $A \in \Rset^{n \times n}$ is said to be 
Metzler if each off-diagonal entry of $A$ is nonnegative. 
For a square matrix $A \in\Rset^{n \times n}$, 
$A^+ = \left(\max\{a_{i,j},0\}\right)_{i,j=1,1}^{n,n}$ is defined, 
where the notation $A = \left(a_{i,j}\right)_{i,j=1,1}^{n,n}$ is used. 
Let $A^- \in\Rset^{n \times n}$ be defined by $A^- = A^+ - A$. 

\section{An Interval Observer}\label{sec:minobs}

This paper considers the following system 
\begin{subequations}\label{eq:sys}
\begin{align}
\dot{x}(t)&= A(y(t))x(t) +\beta(y(t),u(t)) + \delta(t)
\label{eq:sysx}\\
y(t)&=Cx(t) , 
\label{eq:sysy}
\end{align}
\end{subequations}
which is referred to as a plant. 
The vectors $x(t)\in\Rset^n$ and $y(t)\in\Rset^p$ are  
the state and the output, respectively, at time 
$t\in\Rset_+$. The initial condition is denoted by 
$x(0)=x_0$. 
The input $u: \Rset_+\to\Rset^q$ is supposed to belong to $\mathcal{U}$ which 
denotes the set of piecewise continuous functions. 
The disturbance $\delta: \Rset_+\to\Rset^n$ is 
supposed to belong to $\mathcal{D}$ which 
denotes the set of Lebesgue measurable locally essentially bounded 
functions. 
The matrix $C\in\Rset^{p\times n}$ is constant, 
and the functions $A: \Rset^p\to\Rset^{n\times n}$ and 
$\beta: \Rset^p\times\Rset^q\to\Rset^n$ are supposed to be locally 
Lipschitz. 
The maximal open subinterval (of $\Rset_+$) in which the unique $x(t)$ 
exists is denoted by $[0,T_{x_0,u,\delta})$. Hence, 
$T_{x_0,u,\delta}=\sup \{t\in\Rset_+ : |x(t)|<\infty\}$ is 
the escape time for given $x_0$, $u$ and $\delta$. 
System \eqref{eq:sysy} is said to be forward complete 
\cite{ANGSONforcomLyap99} if 
$T_{x_0,u,\delta}=\infty$ holds for any $x_0\in\Rset^n$, 
any $\delta\in\mathcal{D}$ and any $u\in\mathcal{U}$. 
The problem considered in this paper is to design a system 
that generates intervals 
to which individual components of the state $x(t)$ of the plant 
\eqref{eq:sys} belong at each time $t$ based on the information 
of the measurement output 
$y$ and the input $u$. The state $x(t)$ is not measured. 
We assume that the piecewise continuous 
functions $\underline{\delta}, \overline{\delta}^-: \Rset_{\geq 0}\to\Rset^n$ satisfying 
\begin{align}
\underline{\delta}(t)\le\delta(t)\le\overline{\delta}(t)
, \quad \forall t\in\Rset_+   \ \mbox{(a.e.)} 
\label{eq:dist}
\end{align}
are known.

Let $\phi$ is an positive integer which has yet to be determined. 
Define the set $\bfK=\{1,2,\ldots,\phi\}$. 
Let $\overline{L}_k$ and 
$\underline{L}_k$ are 
$\Rset^{n\times p}$ matrices which will be chosen 
for each $k\in\bfK$. 
This paper propose the pair of the differential equations 
\begin{subequations}\label{eq:mintobs}
\begin{align}
&
\dot{\overline{x}}=A(y)\overline{x}
+\overline{Q}(\overline{x},y)
+B(y,u)+\overline{\delta}
\label{eq:mintobsover}
\\
&
\dot{\underline{x}}=A(y)\underline{x}
+\underline{Q}(\underline{x},y)
+B(y,u)+\underline{\delta}  
\label{eq:mintobsunder}
\end{align}
\end{subequations}
as an interval observer, where 
$\overline{Q}$ and $\underline{Q}$ are defined by 
\begin{subequations}\label{eq:mintobsQ}
\begin{align}
&
\overline{Q}_i(\overline{x},y)=\min_{k\in\bfK} [\overline{L}_k]_i(C\overline{x}-y)
\label{eq:mintobsQover}
\\
&
\underline{Q}_i(\underline{x},y)=\max_{k\in\bfK} [\underline{L}_k]_i(C\underline{x}-y)
\label{eq:mintobsQunder}
\end{align}
\end{subequations}
for $i\in\bfI:=\{1,2,\ldots,n\}$. 
Here, the subscripts $i$ denote the $i$-th component 
or the $i$-th row as 
\begin{align}
&
\overline{Q}=\left[\begin{array}{c}
\overline{Q}_1\\
\overline{Q}_2\\[-.4ex]
\vdots\\
\overline{Q}_n
\end{array}\right]
, \quad 
\underline{Q}=\left[\begin{array}{c}
\underline{Q}_1\\
\underline{Q}_2\\[-.4ex]
\vdots\\
\underline{Q}_n
\end{array}\right]
\label{eq:mintobsQpart}
\\
&
\overline{L}_k=
\left[\begin{array}{c}
{[\overline{L}_k]_1}\\ 
{[\overline{L}_k]_2}\\ 
\vdots\\
{[\overline{L}_k]_n}
\end{array}\right]
, \quad 
\underline{L}_k=\left[\begin{array}{c}
{[\underline{L}_k]_1}\\ 
{[\underline{L}_k]_2}\\ 
\vdots\\
{[\underline{L}_k]_n}
\end{array}\right] . 
\label{eq:mintobsLpart}
\end{align}
Due to the minimization taking effect (resp., maximization) 
in \eqref{eq:mintobsQ} for $\phi\neq 1$,  
the function $\overline{Q}$ 
(resp., $\underline{Q}$) becomes nonlinear in $\overline{x}$
(resp., $\underline{x}$). Hence, the observer \eqref{eq:mintobs} 
is nonlinear even when the plant \eqref{eq:sys} is linear. 
The definition \eqref{eq:mintobsQ} also implies that 
the functions $\overline{Q}$ and 
$\underline{Q}$ are Lipschitz in 
$\overline{x}$, $\underline{x}$ and $y$. 
Due to the assumption on $\beta$, for each initial condition, 
the differential equation \eqref{eq:mintobs} 
admits a unique solution $[\overline{x}(t)^T, \underline{x}(t)^T]^T$ 
up to the time when $y$ or $u$ explodes to infinity. 
The following is the main result. 

\begin{thm}\label{thm:maintheo}
Suppose that for each $k\in\bfK$, 
$A(\omega)+L_kC$ is Metzler for all $\omega\in\Rset^p$. 
Assume that there exist integers 
$\overline{l}$, $\underline{l}\in\bfK$, real vectors 
$\overline{v}=[\overline{v}_1,\overline{v}_2,\ldots,\overline{v}_n]^T\gg 0$,
$\underline{v}=[\underline{v}_1,\underline{v}_2,\ldots,\underline{v}_n]^T\gg 0$
and real scalars $\overline{\epsilon}$, $\underline{\epsilon}>0$ such that 
\begin{subequations}\label{eq:righteigv} 
\begin{align}
\forall i\in\bfI\hspace{1.5ex}
\forall \omega\in\Rset^p\hspace{1.5ex} [A(\omega)+\overline{L}_{\overline{l}}C]_i \overline{v}
\le-\overline{\epsilon}\,\overline{v}_i  
\label{eq:righteigvu} 
\\
\forall i\in\bfI\hspace{1.5ex}
\forall \omega\in\Rset^p\hspace{1.5ex} [A(\omega)+\underline{L}_{\underline{l}}C]_i \underline{v}
\le-\underline{\epsilon}\,\underline{v}_i . 
\label{eq:righteigvl} 
\end{align}
\end{subequations}
Then for any $u\in\mathcal{U}$ and 
any $x_0\in\Rset^n$ satisfying 
\begin{align}
\underline{x}(0)\le x_0\le \overline{x}(0) , 
\label{eq:xinit}
\end{align}
system \eqref{eq:mintobs} defined with system \eqref{eq:mintobsQ} 
\eqref{eq:mintobsQpart} and \eqref{eq:mintobsLpart} achieves 
\begin{enumerate}
\renewcommand{\labelenumi}{{\hspace{-1.7ex}\makebox[2.8ex]{\itshape(\roman{enumi})}}}
\renewcommand{\theenumi}{(\roman{enumi})}
\item\label{item:obs} 
If $\delta(t)\equiv\delta^+(t)\equiv\delta^-(t)\equiv0$, 
the implication 
\begin{align}
T_{x_0,u,\delta}\!=\!\infty
\ \Rightarrow\ 
\left\{\begin{array}{l}
\displaystyle\lim_{t\to \infty}|\overline{x}(t) - x(t)| = 0 \\
\displaystyle\lim_{t\to \infty}|\underline{x}(t) - x(t)| = 0 
\end{array}\right\} . 
\label{eq:conv}
\end{align}
\item\label{item:framer} 
If $\delta\in\mathcal{D}$ satisfies \eqref{eq:dist},  
it holds that 
\begin{align}
\underline{x}(t) \leq x(t) \leq \overline{x}(t) 
, \quad \forall t\in [0,T_{x_0,u,\delta}) . 
\label{eq:frame}
\end{align}
\item\label{item:iss}
If $\delta\in\mathcal{D}$ satisfies \eqref{eq:dist}, 
it holds that 
\begin{subequations}\label{eq:issagain}
\begin{align}
&
\limsup_{t\to\infty}\max_{i\in\bfI} 
\frac{1}{\overline{v}_i}(\overline{x}_i(t)-x_i(t))
\nonumber\\
&
\hspace{10ex}\le
\frac{n}{\overline{\epsilon}}\limsup_{t\to\infty}
\max_{i\in\bfI}\frac{\overline{\delta}_i-\delta_i}{\overline{v}_i}  
\label{eq:issagainu}
\\\
&
\limsup_{t\to\infty}\max_{i\in\bfI} 
\frac{1}{\underline{v}_i}(x_i(t)-\underline{x}_i(t))
\nonumber\\
&
\hspace{10ex}\le
\frac{n}{\underline{\epsilon}}\limsup_{t\to\infty}
\max_{i\in\bfI}\frac{\delta_i-\underline{\delta}_i}{\underline{v}_i} . 
\label{eq:issagainl}
\end{align}
\end{subequations}
\end{enumerate}
\end{thm}
Item \ref{item:obs} in Theorem \ref{thm:maintheo} implies that 
each of \eqref{eq:mintobsover} and \eqref{eq:mintobsunder}
has the standard observer property, 
while Item \ref{item:framer} 
is said to be the framer property which is not guaranteed by 
classical standard observers. In fact, 
$\overline{x}$ is the upper frame, while 
$\underline{x}$ is the lower frame. 
A system generating $\overline{x}(t)$ and $\underline{x}(t)$ 
which satisfy \ref{item:obs} and \ref{item:framer} 
is called an interval observer
(see e.g., \cite{DMN-13} and references therein), 
Properties in the form of \eqref{eq:issagainu} and 
\eqref{eq:issagainl} are ultimate 
boundedness \cite{Khalil_Book_02}. They characterize attractivity in the 
presence of disturbances. 
In particular, in the framework of input-to-state stability 
(ISS) \cite{SONISSV}, properties 
\eqref{eq:issagainu} and \eqref{eq:issagainl} are 
called asymptotic gains, which imply that \eqref{eq:mintobs} 
is an ISS interval observer in the definition introduced by 
\cite{ITODINejc18}. 

\begin{rem}
Property \eqref{eq:righteigvu} is a generalized 
expression of requiring $A(\omega)+\overline{L}_{\overline{l}}C$ to 
be Hurwitz uniformly in $\omega$. 
In fact, due to the Perron-Frobenius theorem 
\cite{MeyerMatrixBook,BerPleNonegmat79}, 
in the case where $A$ is constant, 
since $A+\overline{L}_{\overline{l}}C$ is assumed to be Metzler, 
the existence of a vector 
$\overline{v}\gg 0$ and a scalar $\overline{\epsilon}>0$ satisfying \eqref{eq:righteigvu} is 
equivalent to $A+\overline{L}_{\overline{l}}C$ being Hurwitz. 
The same remark applies to $A+\underline{L}_{\underline{l}}C$. 
\end{rem}

\section{The Idea and Proofs}\label{sec:idprf}

\subsection{A Simple Idea: Minimal Error Dynamics}\label{ssec:idea}

The idea of the system 
\eqref{eq:mintobs} proposed with \eqref{eq:mintobsQ} is 
to obtain tighter intervals by making use of 
$\phi$ multiple pairs of the observer gains 
$\overline{L}_k$ and $\underline{L}_k$. 
Define 
\begin{align}
\overline{e}=\overline{x}-x
, \quad 
\underline{e}=x-\underline{x} . 
\label{eq:defe}
\end{align}
Let 
$\underline{e}(t)\!=\!
[\underline{e}_1(t),\underline{e}_2(t),\ldots,\underline{e}_n(t)]^T$ 
and
$\overline{e}(t)\!=\!
[\overline{e}_1(t),\overline{e}_2(t),$ $\ldots,\overline{e}_n(t)]^T$. 
The minimization and the maximization in \eqref{eq:mintobsQ} 
are feasible since $y$ is measured. 
Using the identity $y=Cx$ one can verified from 
\eqref{eq:sys}, \eqref{eq:mintobs} and \eqref{eq:mintobsQ} that 
\begin{subequations}\label{eq:mintobserror}
\begin{align}
&
\dot{\overline{e}}_i=
\min_{k\in\bfK}[A(y)+\overline{L}_kC]_i\overline{e}
, \quad i\in\bfI 
\label{eq:mintobserroru}
\\
&
\dot{\underline{e}}_i=
\min_{k\in\bfK}[A(y)+\underline{L}_kC]_i\underline{e}
, \quad i\in\bfI 
\label{eq:mintobserrorl}
\end{align}
\end{subequations}
in the case of $\delta(t)\equiv\delta^+(t)\equiv\delta^-(t)\equiv0$.
As long as $\overline{e}(t)\ge 0$ holds, the minimization in 
\eqref{eq:mintobserroru} implies that at every moment, 
the upper frame $\overline{x}(t)$ is directed 
to the true value $x(t)$ at the maximum velocity over 
all $\phi$ observer gains $\overline{L}_k$. In the same way, 
as long as $\underline{e}(t)\ge0$ holds,  
equation \eqref {eq:mintobserrorl} implies that at every moment, 
the lower frame $\underline{x}(t)$ is directed 
to the true value $x(t)$ at the maximum velocity over 
all $\phi$ observer gains $\underline{L}_k$. 
If $\phi=1$, i.e., $\bfK=\{1\}$, 
the minimization in \eqref{eq:mintobserror} disappears, 
and the system \eqref{eq:mintobs} reduces to the standard 
Luenberger-type interval observer \cite{GouzeRapaport_EcMo_00,DMN-13}. 
Implementing this idea requires us to confirm 
the non-negativity of $\overline{e}$ and $\underline{e}$, 
the convergence and the ultimate boundedness property of 
$\overline{e}$ and $\underline{e}$ for general $\phi\ge 1$ 
including $\phi=1$ as a special case. 

\subsection{Proof of Item \ref{item:framer} of Theorem \ref{thm:maintheo}}\label{ssec:thmframer}

From \eqref{eq:sys}, \eqref{eq:mintobs} and \eqref{eq:mintobsQ} 
for \eqref{eq:defe} we obtain 
\begin{subequations}\label{eq:mintobserrord}
\begin{align}
&
\dot{\overline{e}}_i=
\min_{k\in\bfK}[A(y)+\overline{L}_kC]_i\overline{e}
+\overline{\delta}-\delta
, \quad i\in\bfI
\label{eq:mintobserrordu}
\\
&
\dot{\underline{e}}_i=
\min_{k\in\bfK}[A(y)+\underline{L}_kC]_i\underline{e}
+\delta-\underline{\delta}
, \quad i\in\bfI . 
\label{eq:mintobserrordl}
\end{align}
\end{subequations}
Since $A(\omega)+\overline{L}_kC$ is Metzler for all 
$\omega\in\Rset^p$, we have 
\begin{align*}
\overline{e}\ge 0 \mbox{ and }
\overline{e}_i=0
\ \Rightarrow\ 
\forall k\in\bfK\hspace{1.5ex}
[A(y)+\overline{L}_kC]_i\overline{e}\ge 0
\end{align*}
for all $i\in\bfI$. 
Due to \eqref{eq:dist}, applying the above property to 
\eqref{eq:mintobserrordu} establishes 
$\overline{e}(t)\ge 0$ for all $t\in\Rset_+$ for any 
$\overline{e}(0)\ge 0$. Thus, the definition \eqref{eq:defe} 
implies the second inequality in \eqref{eq:frame}. 
This same argument for \eqref{eq:mintobserrordl} establishes 
$\underline{e}(t)\ge 0$ and the first inequality in 
\eqref{eq:frame}. 

\subsection{Proof of Item \ref{item:obs} of Theorem \ref{thm:maintheo}}\label{ssec:thmobs}

For each $k\in\bfK$, define 
$\overline{M}_k(y,\overline{e}):
=(A(y)+\overline{L}_kC)\overline{e}$. 
For all $i\in \bfI$, 
let $\overline{k}_i\in\bfK$ denote the integer $k$ 
achieving the minimum in \eqref{eq:mintobsQover}, which is a 
function of $\overline{x}$ and $y$. 
Then \eqref{eq:mintobserroru} is rewritten as 
\begin{align}
\dot{\overline{e}}=
\left[\begin{array}{c}
\overline{M}_{\overline{k}_1}(y,\overline{e}) \\
\overline{M}_{\overline{k}_2}(y,\overline{e}) \\
\vdots\\
\overline{M}_{\overline{k}_n}(y,\overline{e}) 
\end{array}\right]
=:\overline{M}(\overline{x}, y,\overline{e}) 
\label{eq:mintobserroruM}
\end{align}
in the case of $\delta(t)\equiv\delta^+(t)\equiv\delta^-(t)\equiv0$. 
By definition, the function $\overline{M}$ is continuous. 
Let $\overline{l}\in\bfK$, $\overline{v}\gg 0$ and $\overline{\epsilon}>0$ 
be such that \eqref{eq:righteigvu} holds. 
Due to $L_kC\overline{e}=L_k(C\overline{x}-y)$, 
definition \eqref{eq:mintobsQover} implies 
\begin{align}
[A(y)+\overline{L}_{\overline{k}_i}C]_i\overline{e}
\le 
[A(y)+\overline{L}_lC]_i\overline{e} . 
\label{eq:Ldomi}
\end{align}
Define a positive definite and radially unbounded function of 
$\overline{e}\in\Rset_+^n$ as
\begin{align}
\overline{V}(\overline{e})=\max_{i\in\bfI} \frac{1}{\overline{v}_i} \overline{e}_i , 
\label{eq:vmax}
\end{align}
which is in the form of a popular Lyapunov function in monotone systems 
theory (see \cite{DIRR15} and references therein). 
By virtue of \eqref{eq:righteigvu} and \eqref{eq:Ldomi}, 
one can verify 
\begin{align}
\overline{V}(\overline{e})=\frac{1}{\overline{v}_i} \overline{e}_i
\ \Rightarrow \ 
\frac{\partial \overline{V}}{\partial \overline{e}}
\overline{M}(\overline{x},y,\overline{e})
\le-\frac{\overline{\epsilon}}{\overline{v}_i}\overline{e}_i  
\label{eq:vmaxneg}
\end{align}
for \eqref{eq:mintobserroruM}. 
The function $\overline{V}$ defined in \eqref{eq:vmax} is 
locally Lipschitz. 
Let $\calN$ denote the subset of $\Rset^n$ where the gradient 
${\partial \overline{V}}/{\partial \overline{e}}$ does not exist. 
Property \eqref{eq:vmaxneg} guarantees
\begin{align}
\frac{\partial \overline{V}}{\partial \overline{e}}
\overline{M}(\overline{x},y,\overline{e})
\le-\frac{\overline{\epsilon}}{n}\overline{V} , 
\quad \forall \overline{e}\in\Rset_+^n\setminus\calN .
\label{eq:vmaxgas}
\end{align}
Rademacher's theorem tells that the set $\calN$ has measure zero. 
The lower right-hand Dini derivative 
agrees with $({\partial \overline{V}}/{\partial \overline{e}})\overline{M}$ except in 
$\calN$ \cite{BACROSLiapbook05}, i.e., 
the solution $\overline{e}(t)$ of \eqref{eq:mintobserroruM} satisfies 
\begin{align*}
\liminf_{t \rightarrow 0+}
\frac{\overline{V}(\overline{e}(t))-\overline{V}(\overline{e}(0))}{t} 
\le-\frac{\overline{\epsilon}}{n}\overline{V}(\overline{e}(0)) , 
\quad \forall \overline{e}(0)\in\Rset_+^n . 
%\label{eq:vmaxgasdini}
\end{align*}
This proves 
$\lim_{t\to \infty}|\overline{e}(t)| = 
\lim_{t\to \infty}|\overline{x}(t) - x(t)| = 0$ in \eqref{eq:conv}. 
Repeating the same argument for $\overline{x}(t)$ proves 
$\lim_{t\to \infty}|\underline{e}(t)| = 
\lim_{t\to \infty}|x(t) - \underline{x}(t)| = 0$ in \eqref{eq:conv}. 

\subsection{Proof of Item \ref{item:iss} of Theorem \ref{thm:maintheo}}\label{ssec:thmiss}

With $\overline{M}$ in \eqref{eq:mintobserroruM}, 
the error system \eqref{eq:mintobserrordu} is represented by 
\begin{align}
\dot{\overline{e}}=
\overline{M}(\overline{x}, y,\overline{e}) 
+\overline{\delta}-\delta . 
\label{eq:mintobserrorudM}
\end{align}
The function $\overline{V}$ given in \eqref{eq:vmax} satisfies  
\begin{align}
&
\frac{\partial \overline{V}}{\partial \overline{e}}
(\overline{M}(\overline{x},y,\overline{e})
+\overline{\delta}-\delta)
\le-\frac{\overline{\epsilon}}{n}\overline{V}
+\max_{i\in\bfI}\frac{\overline{\delta}_i-\delta_i}{\overline{v}_i}, 
\nonumber\\
&\hspace{31ex}
\forall \overline{e}\in\Rset_+^n\setminus\calN .
\label{eq:vmaxiss}
\end{align}
Therefore, the function $\overline{V}$ is 
an ISS Lyapunov function of system \eqref{eq:mintobserrorudM} 
\cite{SONISSV}, which establishes 
\eqref{eq:issagainu} \cite{SONISSV}. 
The asymptotic gain \eqref{eq:issagainl} of the other error system 
\eqref{eq:mintobserrordu} is obtained in the same say. 

\section{Relaxing the Metzler Assumption}\label{sec:sindef}

The interval observer proposed in this paper 
allows one to use the popular technique of coordinate transformations 
when it is hard or impossible to find $\overline{L}_k$ and 
$\underline{L}_k$ rendering 
$A+\underline{L}_kC$ and $A+\overline{L}_kC$ Metzler. 
The observer can be constructed after applying $z=Rx$ 
to the plant \eqref{eq:sys} for a nonsingular matrix 
$R\in\Rset^{n\times n}$ \cite{RASEFIcoordinate12,EfiRassiSurv16}. 
The upper and lower bounds of disturbances 
are expressed in the new coordinate as 
$R^+\overline{\delta}-R^-\underline{\delta}$ and 
$R^+\underline{\delta}-R^-\overline{\delta}$, respectively. 
In the same way, the initial upper and lower frames
are expressed in the new coordinate as 
$R^+\overline{x}(0)-R^-\underline{x}(0)$ and 
$R^+\underline{x}(0)-R^-\overline{x}(0)$, respectively. 
The upper and lower frames in 
the original coordinate $x$ can be obtained from those 
in the new coordinate $z$ as 
\begin{align}
\overline{x}=S^+\overline{z} - S^-\underline{z}
, \qquad 
\underline{x}=S^+\underline{z} - S^-\overline{z} , 
\label{eq:obsout}
\end{align}
where $S=R^{-1}$. Importantly, the property 
\begin{align}
\underline{x}=\overline{x}  
\ \Leftrightarrow \ 
\underline{z}=\overline{z}
\end{align}
holds. Therefore, the properties achieved in Theorem \ref{thm:maintheo} 
are qualitatively invariant under coordinate changes. 

\section{An Example}\label{sec:exa}

\begin{figure}[t]
%\hspace{-1ex}%
\begin{tabular}{c}
\includegraphics[width=8.2cm,height=4.4cm]{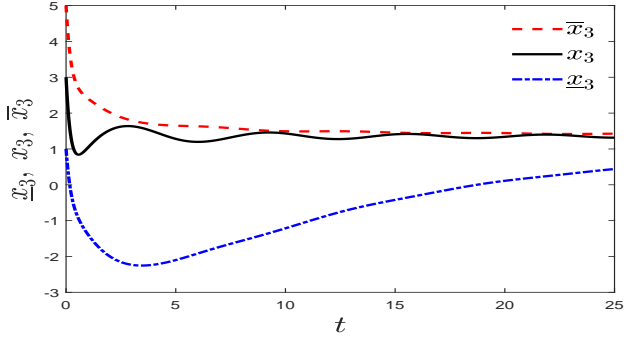}
\\[-1,8ex]
\end{tabular}
\caption{The state $x_3$ and the interval 
$[\underline{x}_3. \overline{x}_3]$ estimated by
the standard interval observer with the single gain $G_1$. 
}
\label{fig:11x3}
%\vspace{2ex}
\end{figure}

\begin{figure}[t]
%\hspace{-1ex}%
\begin{tabular}{c}
\includegraphics[width=8.2cm,height=4.4cm]{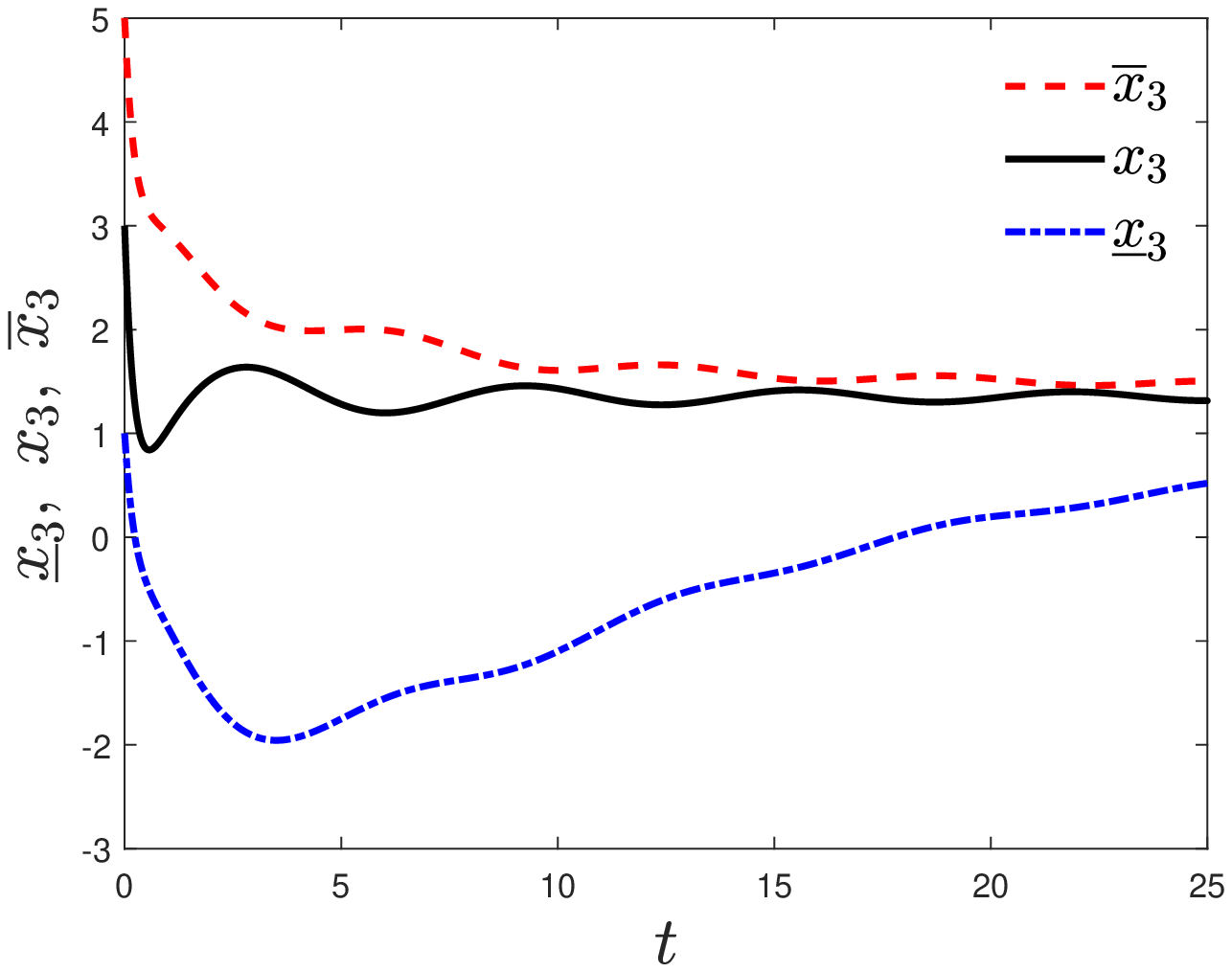}
\\[-1,8ex]
\end{tabular}
\caption{The state $x_3$ and the interval 
$[\underline{x}_3. \overline{x}_3]$ estimated by
the standard interval observer with the single gain $G_2$. 
}
\label{fig:22x3}
%\vspace{2ex}
\end{figure}

\begin{figure}[t]
%\hspace{-1ex}%
\begin{tabular}{c}
\includegraphics[width=8.2cm,height=4.4cm]{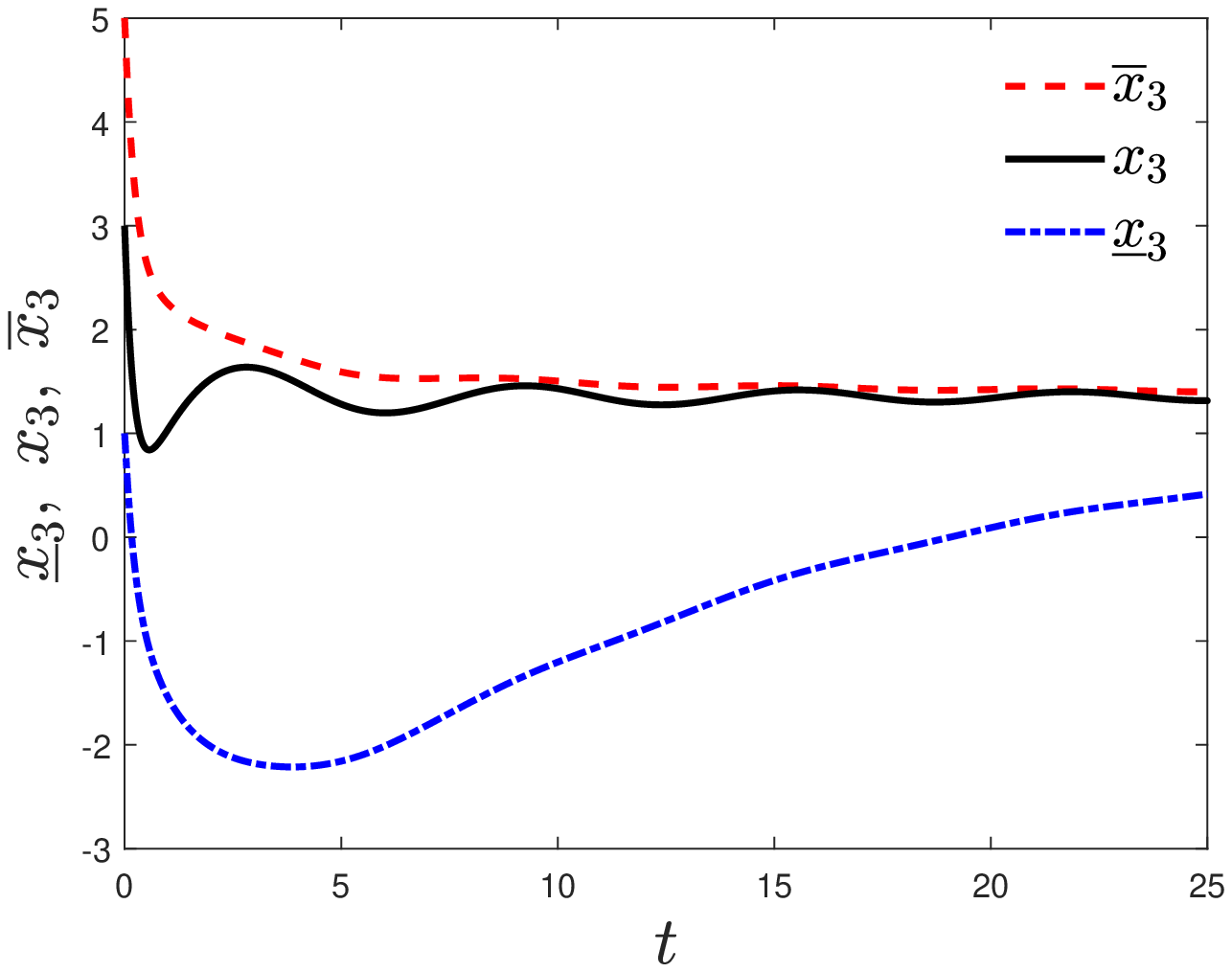}
\\[-1,8ex]
\end{tabular}
\caption{The state $x_3$ and the interval 
$[\underline{x}_3. \overline{x}_3]$ estimated by
the standard interval observer with the single gain $G_3$. 
}
\label{fig:33x3}
%\vspace{2ex}
\end{figure}

\begin{figure}[t]
%\hspace{-1ex}%
\begin{tabular}{c}
\includegraphics[width=8.2cm,height=4.4cm]{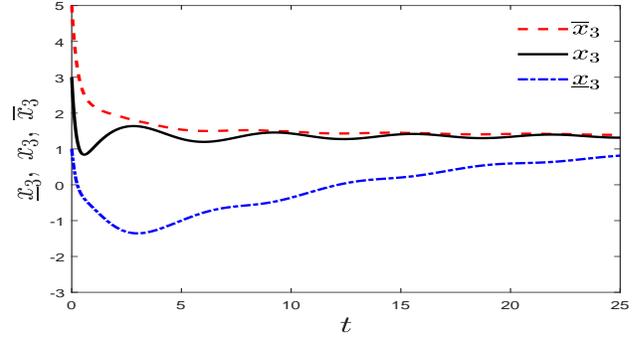}
\\[-1,8ex]
\end{tabular}
\caption{The state $x_3$ and the interval 
$[\underline{x}_3. \overline{x}_3]$ estimated by 
the proposed interval observer 
\eqref{eq:mintobs} with the three gains \eqref{eq:exL123}. 
}
\label{fig:13x3}
%\vspace{2ex}
\end{figure}

Consider the plant \eqref{eq:sys} given by 
\begin{subequations}\label{eq:exsys}
\begin{align}
&
A(y)=\left[\begin{array}{ccc}
-1  & 0.5 & 0   \\[-.5ex]
1   & -1  & 0.8 \\[-.5ex]
0.3 &   1 &  -4
\end{array}\right]
, \quad 
C=\left[\begin{array}{ccc}
1 & 0 & 0 \\
0 & 1 & 0
\end{array}\right] 
\\
&
\beta(y,u) =
\left[\begin{array}{c}
0              \\[-.5ex]
y_2^2-0.2y_2^3 \\[-.5ex]
0
\end{array}\right]
, \quad 
u=0 
. 
\end{align}
\end{subequations}
Define the following gain matrices 
\begin{align*}
&
G_1=\left[\begin{array}{cc}
 -1 &  0 \\[-.5ex]
 0  &  -1 \\[-.5ex]
 -0.3 & -0.3
\end{array}\right], \ 
G_2=\left[\begin{array}{cc}
-0.5 & -0.5 \\[-.5ex]
  -1 &   0 \\[-.5ex]
  0  &  0.2
\end{array}\right], \ 
\\
&
G_3=\left[\begin{array}{cc}
   0 &  -0.5 \\[-.5ex]
   0 &   0   \\[-.5ex]
 0.5 &  -1
\end{array}\right] . 
%\label{eq:exgains}
\end{align*}
Since for each $i=1,2,3$, $A+CG_i$ is Metzler and Hurwitz, 
three interval observers in the standard form 
\cite{GouzeRapaport_EcMo_00,DMN-13}
can be constructed and expressed by \eqref{eq:mintobs} with 
$\phi=1$ in \eqref{eq:mintobsQ}. For 
\begin{align}
&
x_0=\left[\begin{array}{c}
2 \\[-.5ex]
3 \\[-.5ex]
3
\end{array}\right]
, \ 
\overline{x}(0)=\left[\begin{array}{c}
4 \\[-.5ex]
5 \\[-.5ex]
5
\end{array}\right]
, \ 
\underline{x}(0)=\left[\begin{array}{c}
0 \\[-.5ex]
1 \\[-.5ex]
1
\end{array}\right]  
\label{eq:exaset1}
\\
&
\delta(t)=\left[\begin{array}{c}
 \dfrac{2\cos t}{1+t} \\[1ex]
 \dfrac{4\sin t}{1+t} \\[1ex]
-\dfrac{4\cos t}{1+t}
\end{array}\right]  
, \ 
\overline{\delta}(t)=\left[\begin{array}{c}
\dfrac{2}{1+t} \\[1ex]
\dfrac{4}{1+t} \\[1ex]
\dfrac{4}{1+t} 
\end{array}\right]  
, \ 
\underline{\delta}(t)=-\overline{\delta}(t), 
\nonumber\\[-4ex]
&
\label{eq:exaset2}
\end{align}
Figs. \ref{fig:11x3}, \ref{fig:22x3} and \ref{fig:33x3} 
show the upper frame $\overline{x}_3(t)$ and the lower frame 
$\underline{x}_3(t)$ computed by the standard interval 
observer with 
$\overline{L}_1=\underline{L}_1=G_1$, 
$\overline{L}_1=\underline{L}_1=G_2$, 
and 
$\overline{L}_1=\underline{L}_1=G_3$, respectively. 
It is seen that the unmeasured variable $x_3(t)$ of the plant lies 
between the upper and the lower frames. 
The proposed observer \eqref{eq:mintobs} defined with 
\eqref{eq:mintobsQ} and 
\begin{align}
\phi=3, \ 
\overline{L}_1=\underline{L}_1=G_1, \ 
\overline{L}_2=\underline{L}_2=G_2, \ 
\overline{L}_3=\underline{L}_3=G_3 . 
\label{eq:exL123}
\end{align}
fulfills all the assumptions in Theorem \ref{thm:maintheo}. 
For the same initial conditions \eqref{eq:exaset1} and disturbances 
\eqref{eq:exaset2}, 
the upper frame and the lower frame  
generated by the proposed observer 
of the unmeasured variable $x_3(t)$ 
are plotted in Fig. \ref{fig:13x3}. 
The simulation result confirms properties \eqref{eq:frame}, 
\eqref{eq:issagainu} and \eqref{eq:issagainl}. 
The proposed observer gives 
tighter estimates than any of the three single gain 
interval observers. Moreover, the intervals generated by 
the proposed observer are tighter than the intersection of 
the intervals computed at each instant by the three observers.

\section{Conclusions}\label{sec:conc}

The proposed simple observer has clear advantages over the use of the 
intersection of intervals estimated by multiple observers. 
Firstly, the dimension of the observer is fixed and it does not 
increase as one wants to obtain tighter intervals. 
Secondly, the observer does not continue to use loose 
intervals. At every moment, intervals of all state components 
are generated toward the tightest ones. 
The intervals converge to the true points rapidly if disturbances 
vanish. 
These features have been demonstrated by a numerical example. 
The idea is intuitive. This paper has shown how it can be 
justified and implemented.  

%\appendix
\bibliographystyle{plain}
\bibliography{ref}

%=============================================
\end{document}